
\input gtmacros
\input gtmonout
\volumenumber{1}
\volumeyear{1998}
\volumename{The Epstein birthday schrift}
\papernumber{1}
\pagenumbers{1}{21}
\received{10 May 1998}
\published{21 October 1998}

\def\l{[\![}
\def\r{]\!]}
\def\rightangle{\, \vrule height0.25cm width0.4pt depth0.0cm
		   \vrule height0.4pt width0.22cm depth0pt
		   \,}
\def\E{\endprf}

\def\R{\Bbb{R}}
\def\SS{\Bbb{S}}
\def\MM{{\cal M}}
\def\NN{{\cal N}}
\def\TT{{\cal T}}
\let\qua\stdspace

\reflist

\refkey\AW {\bf W\,K Allard}, {\it On the first variation of a varifold},
       Annals of Math. {95} (1972) 417--491

\refkey\BG {\bf M~Berger}, {\bf B~Gostiaux}, {\it Geometrie 
differentielle: varietes, courbes et surfaces}, 
Presses Universitaires de France, III (1987)

\refkey\T {\bf J~Cheeger}, {\bf W~M\"uller}, {\bf R~Schrader}, {\it On the 
curvature of
piecewise flat spaces}, Comm.~Math.~Phys. {92} (1984) 405--454

\refkey\CSW {\bf R Connelly}, {\bf I Sabitov}, {\bf A Walz},
{\it The bellows conjecture},
Beitr\"age Algebra Geom. 38 (1997) 1--10 

\refkey\H {\bf G~Herglotz}, {\it Ueber der Starrheit der Eiflachen}, 
Abh. Math. Semin. Hansische Univ. {92} (1943) 127--129

\refkey\RS {\bf I~Rivin}, {\bf J-M~Schlenker}, {\it Schl\"afli formula 
and Einstein
manifolds}, IHES pre\-print (1998)

\refkey\S 
{\bf D\,V\, Alekseevsky}, {\bf \`E\,B Vinberg}, {\bf A\,S
Solodovnikov}, {\it Geometry of spaces of constant curvature}, from:
``Geometry II'', Encyclopaedia Math. Sci. 29, Springer--Verlag,
Berlin (1993)

\refkey\SM {\bf M Spivak}, {\it A Comprehensive Introduction to Differential
Geometry}, (Second Edition) Publish or Perish, Berkeley (1979)

\endreflist

\title{The mean curvature integral is invariant\\under bending}

\author{Frederic J Almgren Jr\\Igor Rivin}

\address{Mathematics Institute, University of Warwick\\
Coventry, CV4 7AL, UK}
\email{igor@maths.warwick.ac.uk}

\abstract
Suppose $\MM_t$ is a smooth family of compact connected two
dimensional submanifolds of Euclidean space $E^3$ without boundary varying
isometrically in their induced Riemannian metrics.  Then we show that the
mean curvature integrals\vglue -4mm
$$\int_{\MM_t}H_t\,d {\cal H}^2$$
are constant.  It is unknown whether 
there are nontrivial such bendings $\MM_t$.  The estimates
also hold for periodic manifolds for which there are nontrivial
bendings. In addition, our methods work essentially without change to
show the similar results for submanifolds of $H^n$ and $S^n$, to wit,
if $\MM_t = \partial X_t$\vglue -4mm
$$d \int_{\MM_t}H_t\,d {\cal H}^2 = -k {n-1} d V(X_t),$$
where $k=-1$ for $H^3$ and $k=1$ for $S^3$. The Euclidean case can be
viewed as a special case where $k=0$.
The rigidity of the mean curvature integral can be used to show new
rigidity results for isometric embeddings and provide new proofs of
some well-known results. This, together with far-reaching extensions
of the results of the present note is done in the preprint [\RS]. Our
result should be compared with the well-known formula of Herglotz (see
[\H], also [\SM] and [\BG]). 
\endabstract

\asciiabstract{%
Suppose M_t is a smooth family of compact connected two dimensional
submanifolds of Euclidean space E^3 without boundary varying
isometrically in their induced Riemannian metrics.  Then we show that
the mean curvature integrals over M_t are constant.  It is unknown
whether there are nontrivial such bendings.  The estimates also hold
for periodic manifolds for which there are nontrivial bendings. In
addition, our methods work essentially without change to show the
similar results for submanifolds of H^n and S^n.  The rigidity of the
mean curvature integral can be used to show new rigidity results for
isometric embeddings and provide new proofs of some well-known
results. This, together with far-reaching extensions of the results of
the present note is done in the preprint: I Rivin, J-M Schlenker,
Schlafli formula and Einstein manifolds, IHES preprint (1998). Our
result should be compared with the well-known formula of Herglotz.}




\maketitle

\section{Introduction}

The underlying idea of this note is the following.  Suppose
$\NN_t$ is a smoothly varying family of polyhedral solids
having edges $\big\{E_t(k)\big\}_k$, and associated
(signed) dihedral angles $\big\{\theta_t(k)\big\}_k$.  According to
a theorem of Schlafli [\S]
$$\sum_{k}\big|E_t(k)\big|\,
{d\over dt}\theta_t(k) = 0.$$
In case edge length is preserved in the family, ie
$${d \over dt}\big| E_t(k)\big| = 0$$
for each time $t$ and each $k$, then also (product rule)
$${d \over dt}\sum_{k} \big| E_t(k)\big|
\,\theta_t(k) = 0.$$
Should the $\partial \NN_t$'s be polyhedral approximations to
submanifolds $\MM_t$ varying isometrically, one might regard
$$\sum_{k}\big|E_t(k)\big|\, \theta_t(k)$$ 
as a reasonable approximation to the mean curvature integrals
$$\int_{\MM_t} H_t\, d {\cal H}^2$$
and expect
$${d \over dt}\big| E_t(k)\big|$$
to be small.  Hence it is plausible that the mean curvature
integrals of the $\MM_t$'s might be constant.
In this note we show that that is indeed the case.  

Examples such as the isometry pictured on page 306 of volume 5
of [\SM] show that the mean curvature integral is not preserved
under discrete isometries.

Two comments are in order. The first is that it is very likely that
there are {\it no} isometric bendings of hypersurfaces. One reason for
the existence of the current work is to produce a tool for resolving
this conjecture (as Herglotz' mean curvature variation formula can be
used to give a simple proof of Cohn--Vossen's theorem on rigidity of
convex hypersurfaces). Secondly, the main theorem can be viewed as a
sort of dual bellows theorem (when the hypersurface in question lies
in $H^n$ or $S^n$): as the surface is isometrically deformed, the
volume of the {\sl polar dual} stays constant. This should be
contrasted with the usual bellows theorem recently proved by
Sabitov,  Connelly and Walz [\CSW].

\section{Terminology and basic facts}

Our object in this section is to set up terminology for a family of
manifolds varying smoothly through isometries.  We consider triangulations
of increasing fineness varying with the manifolds.  To make possible
our mean curvature analysis we associate integral varifolds with both
the manifolds and the polyhedral surfaces determined by the
triangulations.  The mean curvature integral of interest is identified
with (minus two times)
the varifold first variation associated with the unit normal
initial velocity vector field.  

\sh{2.1\qua {T}erminology and {f}acts
for a {s}tatic {m}anifold {$\MM$}}

{\bf 2.1.1}\qua  We suppose that $\MM\subset \R^3$ is a compact connected 
smooth two dimensional 
submanifold of $\R^3$ without boundary oriented by a smooth Gauss
mapping ${\bf n}\co \MM \to \SS^2$ of unit normal vectors.

\medskip
{\bf 2.1.2}\qua   $H \co \MM \to \R$ denotes half the sum of principal 
curvatures in direction ${\bf n}$ at points in $\MM$
so that $H{\bf n}$ is the mean curvature vector field of $\MM$.

\medskip
{\bf 2.1.3}\qua   We denote by $U$ a suitable neighborhood of $\MM$ in $\R^3$
in which a smooth nearest point retraction mapping
$\rho \co U \to \MM$ is well defined.
The smooth signed distance function $\sigma \co U \to \R$ is defined
by requiring $p = \rho(p) + \sigma(p)\, {\bf n}(\rho(p))$ for each $p$.
We set
 $$g = \nabla \sigma \co U \to \R^3$$
(so that $g|\MM = {\bf n}$); the vector field $g$ is the initial velocity
vector field of the deformation
$$G_t\co U\to \R^3, \quad G_t(p) = p + t \,g(p) \quad\hbox{for } p \in
U.$$

\medskip
{\bf 2.1.4}\qua   We denote by
 $$V = {\bf v}(\MM)$$
the {\sl integral varifold} associated with $\MM$
[\AW, 3.5].  The first variation
distribution of $V$ [\AW, 4.1, 4.2] is representable by integration 
[\AW, 4.3] and can be written
$$\delta V = {\cal H}^2 \rightangle \MM \wedge (-2 H){\bf n}$$
[\AW, 4.3.5]
so that
$$\delta V(g) = {d \over dt}{\cal H}^2\big(G_t(\MM)\big)\bigg|_{t=0}
= -2\int_\MM g \cdot  H \,{\bf n} \, d {\cal H}^2
= -2\int_\MM H \, d {\cal H}^2;$$
here ${\cal H}^2$ denotes two dimensional Hausdorff measure in $\R^3$.

\medskip
{\bf 2.1.5}\qua  By a {\bf vertex} $p$ in $\MM$ we mean any point $p$ in $\MM$.
By an {\bf edge} $\langle pq \rangle$ in $\MM$ we mean any (unordered)
pair of distinct vertexes $p$, $q$ in $\MM$ which are close 
enough together 
that there is a unique length minimizing geodesic arc $\l pq\r$
in $\MM$ joining them; in particular $\langle pq\rangle = \langle qp\rangle$.  
For each edge $\langle pq\rangle$ we write $\partial \langle pq\rangle = \{ p, \, q\}$ 
and call $p$ a vertex of edge $\langle pq\rangle$, etc.  
We also denote
by $\overline{pq}$ the straight line segment in $\R^3$ between
$p$ and $q$, ie the convex hull of $p$ and $q$.
By a {\sl facet} $\langle pqr\rangle$ in $\MM$ we mean any (unordered) triple
of distinct vertexes $p$, $q$, $r$ which are not collinear in $\R^3$
such that $\langle pq\rangle$, $\langle qr\rangle$, $\langle rp\rangle$ are edges in $\MM$; in particular,
$\langle pqr\rangle = \langle qpr\rangle = \langle rpq\rangle $, etc.  For each facet $\langle pqr\rangle$ we write 
$\partial \langle pqr\rangle =\big\{ \langle pq\rangle, \, \langle qr\rangle,\, \langle rp\rangle\big\}$ and call
$\langle pq\rangle$ an edge of facet $\langle pqr\rangle$
and also denote by
$\overline{pqr}$ the convex hull of $p$, $q$, $r$ in $\R^3$.

\medskip
{\bf 2.1.6}\qua  Suppose $ 0 < \tau < 1$ and $ 0 < \lambda < 1 $.
By a $\tau,\lambda$ {\sl regular triangulation $\TT$ of $\MM$ of 
maximum edge length} $L$ we mean

(i)\qua  a family $\TT_2$ of facets in $\MM$, together with

(ii)\qua  the family $\TT_1$ of all edges of facets in $\TT_2$
together with

(iii)\qua  the family $\TT_0$ of all vertexes of edges in $\TT_1$

such that

(iv)\qua  $\overline{pqr}\subset U$ for each facet $\langle pqr\rangle$ in 
$\TT_2$

(v)\qua  $\MM$ is partitioned by the family of subsets
$$\eqalign{\bigg\{\rho\big(\overline{pqr} \sim (\overline{pq} \cup \overline{qr}
\cup \overline{rq})\big):
\langle  pqr\rangle \in \TT_2\bigg\} &\cup
\bigg\{\rho(\overline{pq}) \sim \{p, \, q\} :
\langle  pq\rangle\in \TT_1\bigg\}\cr
 &\cup
\bigg\{\{p\} : p \in \TT_0\bigg\}\cr}$$

(vi)\qua  for facets $\langle pqr\rangle \in \TT_2$ we have the uniform
nondegeneracy
condition: if we set $u = q-p$ and $v = r-p$ then
$$\left|v - \left({u \over |u| }\cdot v\right) {u \over |u|}\right|
\ge \tau|v|$$

(vii)\qua   $L = \sup\big\{|p-q|: \langle pq\rangle\in \TT_1\big\}$

(viii)\qua  for edges in $\TT_1$ we have the uniform control on
the ratio of lengths: 
$$\inf\big\{|p-q|: \langle pq\rangle\in \TT_1\big\}
\ge \lambda L.$$

\medskip
{\bf 2.1.7
{F}act}\qua [\T]\qua It is a standard fact about the geometry of smooth
submanifolds that there are $0 < \tau < 1$ and $0 < \lambda < 1$ such
that for arbitrarily small maximum
edge lengths $L$ there are $\tau, \lambda$ regular 
triangulations of $\MM$ of maximum edge length $L$.
We fix such $\tau$ and $\lambda$.  We hereafter consider only $\tau,
\lambda$ regular triangulations $\TT$ 
 with very small maximum edge length $L$.
Once $L$ is small the triangles $\overline{pqr}$ associated with
$\langle pqr \rangle$ in $\TT_2$ are very nearly parallel with the
tangent plane to $\MM$ at $p$.

\medskip
{\bf 2.1.8}\qua  
Associated with each facet $\langle pqr\rangle$ in $\TT_2$ is the
{\sl unit normal vector} ${\bf n}(pqr)$ to $\overline{pqr}$ having 
positive inner product with the normal ${\bf n}(p)$ to $\MM$ at $p$.

\medskip
{\bf 2.1.9}\qua  Associated with each edge $\langle pq\rangle$ in 
$\TT_1$ are exactly two distinct facets 
$\langle pqr\rangle$ and $\langle pqs\rangle$ in $\TT_2$.  We denote by
$${\bf n}( pq)= 
{{\bf n}( pqr) + {\bf n}( pqs)\over
\big|{\bf n}\big(pqr)+{\bf n}( pqs)\big|} $$
the {\sl average normal vector} at $\overline{pq}$.  

For each $\langle pq\rangle$ we further denote by $\theta(pq)$ 
the {\sl signed dihedral angle}
at $\overline{pq}$ between the oriented plane directions of 
$\overline{pqr}$ and $\overline{pqs}$ which is
characterized by the condition
$$2 \sin \left({\theta(pq)\over 2}\right)\, {\bf n}(pq) =V + W$$
where

$\bullet$\qua  
$V$ is the unit exterior normal vector to $\overline{pqr}$ along edge
$\overline{pq}$, so that, in particular,
$$V \cdot (p-q) = V \cdot {\bf n}(pqr) = 0;$$
$\bullet$\qua 
$W$ is the unit exterior normal vector to $\overline{pqs}$ along edge
$\overline{pq}$.

One checks that
$$\cos \theta(pq) = {\bf n}(pqr) \cdot {\bf n}(pqs).$$
Finally for each $\langle pq \rangle$ we denote by
$$g(pq) = |p-q|^{-1} \int_{\overline{pq}} g\, d {\cal H}^1\in\R^3 $$
the $\overline{pq}$ {\sl average of} $g$; here ${\cal H}^1$ is one dimensional
Hausdorff measure in $\R^3$.

\medskip
{\bf 2.1.10}\qua Associated with our triangulation $\TT$ of $\MM$
is the {\sl polyhedral approximation} 
$$\NN[\TT] = \cup\big\{\overline{pqr} : \langle pqr\rangle \in \TT_2\big\}$$
and the integral varifold
$$V[\TT]=\sum_{\langle pqr \rangle \in \TT_2}
{\bf v}\big(\overline{pqr}\big)
= {\bf v}\big(\NN(\TT)\big)$$
whose first variation distribution is representable by integration
$$\delta V[\TT] = \sum_{\langle pq \rangle \in \TT_1}
{\cal H}^1 \rightangle \overline{pq} \wedge\left[ 2\,
 \sin \left({\theta(pq)\over 2}\right)\right]\,{\bf n}(pq)$$
[\AW, 4.3.5] so that
$$\delta V[\TT](g) = \sum_{\langle pq \rangle \in \TT_1}
\bigg[ |p-q|\bigg]\,\bigg[2\, \sin \left({\theta(pq)\over 2}\right)\bigg]
\, \bigg[{\bf n}(pq)\cdot g(pq)\bigg].$$

\sh{{2.2\qua T}erminology and facts for a {f}low of {m}anifolds {$\MM_t$}}

\medskip
{\bf 2.2.1}\qua  As in 2.1.1 we suppose that
$\MM\subset \R^3$ is a compact connected smooth two dimensional 
submanifold of $\R^3$ without boundary oriented by a smooth Gauss
mapping ${\bf n}\co \MM \to \SS^2$ of unit normal vectors.
We suppose additionally
that $\varphi\co(-1, \, 1)\times \MM \to \R^3$ is 
a smooth mapping with $\varphi(0, \, p) = p$
for each $p \in \MM$.  For each $t$ we set 
$$\varphi[t] = \varphi(t,\,\cdot)\co\MM \to \R^3\quad\hbox{ and }\quad
\MM_t = \varphi[t](\MM).$$  Our principal assumption is that,
for each $t$, the mapping
$\varphi[t]\co \MM \to \MM_t$ {\sl is an orientation preserving
isometric imbedding (of Riemannian manifolds).}
In particular, each
$\MM_t\subset \R^3$ is a compact connected smooth two dimensional 
submanifold of $\R^3$ without boundary oriented by a smooth Gauss
mapping ${\bf n_t}\co \MM_t \to \SS^2$ of unit normal vectors.

\medskip
{\bf 2.2.2}\qua  As in 2.1.2, for each $t$, we denote by $H_t{\bf n}_t$ the 
mean curvature vector field of $\MM_t$.  

{\bf 2.2.3}\qua  As in 2.1.3, for each $t$ we denote by $U_t$ a suitable 
neighborhood of $\MM_t$ in $\R^3$
in which a smooth nearest point retraction mapping
$\rho_t \co U_t \to \MM_t$ is well defined together with
smooth signed distance function $\sigma_t \co U_t \to \R$;
also we set $g[t] = \nabla \sigma_t \co U_t \to \R^3$ as an initial
velocity vector field.

\medskip
{\bf 2.2.4}\qua  By a convenient abuse of notation we assume that we can
define a smooth map
$$\varphi\co (-1, \, 1)\times U_0 \to \R^3,$$
$$\varphi(t, \, p) = 
\varphi\big(t, \,\rho_0(p) + \sigma_0(p) {\bf n}_0(\rho(p)\big)
=\varphi\big(t, \,\rho_0(p)\big) + \sigma_0(p) {\bf n}_t(\rho_0(p)\big)$$
for each $t$ and $p$.
With $\varphi[t] = \varphi(t, \, \cdot)$ we have
$\varphi[0] = {\bf 1}_{U_0}$ and, additionally,
$\sigma_0(p) = \sigma_t\big(\varphi[t](p)\big)$.
We further assume that 
$$U_t = \varphi[t]\,U_0$$
for each~$t$.

\medskip
{\bf 2.2.5 Fact}\qua
If we replace our initial $\varphi[t]\co \MM\to \R^3$'s 
by $\varphi[\mu t]$ for large enough $\mu$
(equivalently, restrict times $t$ to $-1/\mu < t < 1/\mu$)
and decrease the size of $U_0$ then the extended 
$\varphi[t]\co U_0 \to \R^3$'s will exist. 
Such restrictions
do not matter in the proof of our main assertion, since it is local
in time and requires only small neighborhoods of the $\MM_t$'s.

\medskip
{\bf 2.1.6}\qua  As in 2.1.4, for each $t$ we denote by
$$V_t = {\bf v}(\MM_t)$$
the integral varifold associated with $\MM_t$.

\medskip
{\bf 2.2.7}\qua  We fix $0<\tau<1/2$ and $0<\lambda<1/2$ as in 2.1.7 and fix
$2\tau$, $2\lambda$ regular triangulations $\TT(1)$, $\TT(2)$, $\TT(3)$,
$\ldots$ of $\MM$ having maximum edge lengths 
$L(1)$, $L(2)$, $L(3)$ $\ldots$ respectively with 
$\lim_{j \to \infty}L(j) = 0.$
For each $j$, the  vertexes of $\TT(j)$ are denoted $\TT_0(j)$,
the edges are denoted $\TT_1(j)$, and
the facets are denoted $\TT_2(j)$.
For all large $j$ and each $t$ we have
triangulations $\TT(1,\,t)$, $\TT(2,\,t)$, $\TT(3,\,t)$, $\ldots$ 
of $\MM_t$ as follows.
With notation similar to that above we specify, for each $j$ and $t$, 
$$\TT_0(j, \, t) = \bigg\{\varphi[t](p): p \in \TT_0(j)\bigg\},
\quad
\TT_1(j, \, t)=\bigg\{\big\langle\varphi[t](p)\,\varphi[t](q)\big\rangle 
: \langle p q\rangle\in \TT_1(j)\bigg\},$$
$$\TT_2(j, \, t)=\bigg\{\big\langle\varphi[t](p)\, \varphi[t](q)\, 
\varphi[t](r)\big\rangle : \langle pqr\rangle\in \TT_2(j)\bigg\}.$$

\medskip
{\bf 2.2.8 {F}act}\qua
If we replace $\varphi[t]$ by $\varphi[\mu t]$ for large enough $\mu$
(equivalently, restrict times $t$ to $-1/\mu < t < 1/\mu$) then 
$\TT(1,\, t)$, $\TT(2,\, t)$, $\TT(3,\, t)$, $\ldots$ 
will a sequence of $\tau, \lambda$ regular triangulations  
of $\MM$ with maximum edge lengths $L(j, \, t)$ converging to $0$ 
uniformly in time $t$ as $j \to \infty$.  Such restrictions
do not matter in the proof of our main assertion, since it is local
in time.  We assume this has been done, if necessary, and that each
of the triangulations $\TT(j, \, t)$ is
$\tau, \lambda$ regular with maximum edge lengths $L(j, \,t)$
converging to $0$ as indicated.

\medskip
{\bf 2.2.9}\qua  
As in 2.1.8 we associate with each $j$, $t$, and
$\langle pqr\rangle  \in \TT_2(j)$ a unit normal vector
${\bf n}[t, \, j](pqr)\,$ to $\,\overline{\varphi[t](p)\,\varphi[t](q)\,
\varphi[t](r)}\,$.
As in 2.1.9 we associate with each $j$, $t$, and
$\langle pq\rangle  \in \TT_1(j)$ an average normal vector
$\,{\bf n}[t,\,j](pq)\,$ at 
$\,\overline{\varphi[t](p)\,\varphi[t](q)}\,$
and a signed dihedral angle
$\,\theta[t,\,j](pq)\,$ at 
$\,\overline{\varphi[t](p)\,\varphi[t](q)}\,$
and the $\,\overline{\varphi[t](p)\,\varphi[t](q)}\,$
average $\,g[t, \, j](pq)\,$ of $g[t]$.
\medskip
{\bf 2.2.10}\qua   As in 2.1.10 we associate with each triangulation 
$\TT(j, \, t)$ of $\MM_t$ a {\sl polyhedral approximation} 
$\NN[\TT(j,\,t)]$ and an integral varifold
$$V[\TT(j,\,t)] = {\bf v}\big(\NN[\TT(j, \, t)]\big)
= \sum_{\langle pqr \rangle \in \TT_1(j)}
{\bf v}\bigg(\overline{\varphi[t](p)\,\varphi[t](q)\,\varphi[t](r)}\bigg)$$
with first variation distribution 
$$\delta V[\TT(j, \, t)] = \sum_{\langle pq \rangle \in \TT_1(j)}
{\cal H}^1 \rightangle 
\bigg[\overline{\varphi[t])p)\,\varphi[t](q)}\bigg]
\wedge\bigg[ 2\, \sin \left({ \theta[t,\,j](pq) 
\over 2}\right)\bigg]\,{\bf n}[t, \, j](pq).$$
so that
$$\eqalign{\delta &V[\TT(j, \, t)]\big(g[t]\big)\cr 
= &\sum_{\langle pq \rangle \in \TT_1(j)}
\bigg[ \big|\varphi[t](p) - \varphi[t](q)\big|
\bigg]\,\bigg[2\, \sin 
\left({\theta[t,\,j](pq)\over 2}\right)\bigg]
\, \bigg[{\bf n}[t, \, j](pq)\cdot g[t, \, j](pq)\bigg].\cr}$$
 
\medskip
{\bf 2.2.11}\qua   The quantity we wish to show is constant in time is 
$$\int_{\MM_t} H_t\, d {\cal H}^2 = -\left({1\over 2}\right)
\delta V_t\big(g[t]\big).$$
Since, for each time $t$, 
$$V_t = \lim_{j \to \infty}V[\TT(j, \, t)]\qquad\hbox{(as varifolds)}$$
we know, for each $t$,
$$\delta V_t\big(g[t]) = \lim_{j \to \infty}
\delta V[\TT(j, \, t)]\big(g[t]\big).$$
We are thus led to seek to estimate
$${d \over dt} \delta V[\TT(j, \, t)]\big(g[t]\big)$$
using the formula in 2.2.10.  A key equality it provided by
Schlafli's theorem mentioned above which, in the present terminology,
asserts for each $j$ and $t$,
$$\sum_{\langle pq \rangle \in \TT_1(j)}
\bigg[ \big|\varphi[t](p) - \varphi[t](q)\big| \bigg]
\,{d\over dt} \bigg[ \theta[t,\,j](pq)\bigg]=0.$$

\medskip
{\bf 2.2.12 {F}act}\qua   Since, for each $\langle ppq \rangle$
in $\TT_2(j)$,  $\partial \langle pqr \rangle$ consists of exactly
three edges, and, for each $\langle pq \rangle$ in $\TT_1(j)$, 
there are exactly two distinct facets $\langle pqr \rangle$ in 
$\TT_2(j)$ for which $\langle pq \rangle \in \partial \langle 
pqr \rangle$ we infer that, for each $j$, 
$$card \big[\TT_1(j)\big] =  {3\over 2}card \big[\TT_2(j)\big].$$
We then use the $\tau, \lambda$ regularity of the the $\TT(j)$'s
to check that that, for each time $t$ and each
$\langle ppq \rangle$ in $\TT_2(j)$ the following four numbers have
bounded ratios (independent of $j$, $t$, and $\langle ppq \rangle$)
with each other
$${\cal H}^2\bigg(\overline{\varphi[t](p)\,\varphi[t](q)\, \varphi[t](r)}
\bigg), \qquad
\big|\varphi[t](p)-\varphi[t](q)\big|^2, \qquad
L(j, \, t)^2,\qquad L(j)^2.$$
Since
$$\lim_{j \to \infty}{\cal H}^2\big(\NN[j, \, t]\big) =
{\cal H}^2\big(\MM_t\big) ={\cal H}^2\big(\MM\big),$$
we infer
$$\sup_{j}\sum_{\langle pq \rangle \in \TT_1(j)}
L(j)^2 < \infty, \qquad
\lim_{j \to \infty} \sum_{\langle pq \rangle \in \TT_1(j)}
L(j)^3 = 0.$$

\section{Modifications of the flow}

\sh{{3.1\qua  J}ustification for {c}omputing with {m}odified
 {f}lows}

As indicated in 2.2, we wish to estimate the time derivatives of
$$\eqalign{\delta &V[\TT(j, \, t)]\big(g[t]\big)\cr 
=& \sum_{\langle pq \rangle \in \TT_1(j)}
\bigg[ \big|\varphi[t](p) - \varphi[t](q)\big|
\bigg]\,\bigg[2\, \sin \left({\theta[t,\,j](pq)\over 2}\right)\bigg]
\, \bigg[{\bf n}[t, \, j](pq)\cdot g[t, \, j](pq)\bigg].\cr}$$
In each of the $\langle pq \rangle$ summands, each
of the three factors
$$\bigg[ \big|\varphi[t](p) - \varphi[t](q)\big|
\bigg],\quad\bigg[2\, \sin \left({\theta[t,\,j](pq)\over 2}\right)\bigg]
,\quad \bigg[{\bf n}[t, \, j](pq)\cdot g[t, \, j](pq)\bigg]$$
is an intrinsic geometric quantity (at each time) whose value does not
change  under isometries of the ambient $\R^3$.  With
$\langle pqr \rangle$ and $\langle pqs \rangle$ denoting the two facets
sharing edge $\langle pq \rangle$, we infer that each of the factors
depends at most on the relative positions 
of $\varphi[t](p)$, $\varphi[t](q)$, $\varphi[t](r)$,
$\varphi[t](s)$ and $\varphi[t]\MM$.
Suppose $\psi \co(-1, \, 1) \times \R^3 \to \R^3$ is 
continuously differentiable, and for each $t$, the function
$\psi[t] = \psi(t, \,\cdot)\co \R^3 \to \R^3$ is an isometry.
Suppose further, we set
$$\varphi^\ast(t, \, p) = \psi\big(t,\, \varphi(t, \, p)\big), \quad
\varphi^\ast[t]=\varphi^\ast(t, \, \cdot)$$
for each $t$ and $p$ so that $\varphi^\ast[t] 
= \psi[t]\circ\varphi[t]$.  
If we replace $\MM$ by $\MM^\ast= \psi[0]\MM$ and $\varphi$ by
$\varphi^\ast$ then we could follow the procedures of 2.1 and 2.2
to construct triangulations and polyhedral approximations 
$\TT^\ast[j, \,t]$ and varifolds $V^\ast$, etc.  
with
$$ \delta V[\TT(j, \, t)]\big(g[t]\big)=
\delta V^\ast[\TT^\ast(j, \, t)]\big(g^\ast[t]\big).$$
Not only do we have equality in the sum, but, for each $\langle 
pq\rangle$ the corresponding summands are identical numerically.
Hence, in evaluating
$ \delta V[\TT(j, \, t)]\big(g[t]\big)$ we are free to (and will)
use a different $\psi$ and $\varphi^\ast$ for each summand.

\sh{{3.2\qua C}onventions for {d}erivatives}

Suppose $W$ is an open subset of $\R^M$ and 
$f = \big(f^1, \, f^2, \, \ldots, \, f^N\big) \co W \to \R^N$ 
is $K$ times continuously differentiable.  We denote by
$$|||D^{K}f|||$$
the supremum of the partial derivatives 
$${\partial^k f^K\over 
\partial x_{i(1)} \partial x_{i(2)}\ldots \partial x_{i(K)}}(p)$$ 
corresponding to all points $p\in W$, all 
$\big\{ i(1), \, i(2), \, \ldots , \,  i(K)\big\}\subset
\big\{ 1, \, \ldots, \, M\big\}$ and
$k = 1, \, \ldots, N$,
all choices of orthonormal coordinates 
$(x_1, \, \ldots, \, x_M)$ 
for $\R^M$ and all choices of orthonormal coordinates 
$(y_1, \, \ldots, \, y_N)$
for $\R^N$.

\sh{{3.3\qua C}onventions for {i}nequalities}

In making various estimates we will use use the largest edge
length of the $j$th triangulation, typically called $L$, and
a general purpose constant $C$.  The constant $C$ will have
different values in different contexts (even in the same
formula).  What is implied is that, with $\MM$ and $\varphi$
fixed, the constants $C$ can be chosen independent of the
level of triangulation (once it is fine enough) and independent
of time $t$ and independent of the various modifications of our
flow which are used in obtaining our estimates.
As a representative example of our terminology, the expression
$$ A = B \pm C  L^2$$
means 
$$-CL^2 \le A-B \le C L^2.$$

\sh{{3.4\qua F}ixing a {v}ertex at the {o}rigin}

Suppose $p$ is a vertex in $\MM$ and 
$$\varphi_\ast (-1, \, 1) \times U_0 \to \R^3, \qquad
\varphi_\ast(t, \, q) = \varphi(t, \, q) - \varphi(t, \, p)\quad
\hbox{ for each } q.$$
Then $\varphi^\ast(t, \, p) = (0, \, 0, \, 0)$ for each $t$.
One checks, for $K = 0, \, 1, \, 2, \, 3$ that
$$|||\,D^{K}\varphi^\ast\,|||\le 2 |||\,D^{K}\varphi\,|||,\qquad
|||\,D^{K}\varphi^\ast[t]\,||| = |||\,D^{K} \varphi[t]\,\,|||$$
for each $t$.

\sh{{3.5\qua  M}apping a {f}rame to the {b}asis {v}ectors}

Suppose $(0, \, 0, \, 0)\in \MM$ and that ${\bf e}_1$ and
${\bf e}_2$ are tangent to $\MM$ at $(0, \, 0, \, 0)$.
Suppose also
$\varphi(t, \,0, \, 0, \, 0) = (0, \, 0, \, 0)$ for each $t$.
Then the mapping $\varphi^\ast$ given by setting
$$\varphi^\ast[t] 
=\pmatrix{{\partial \varphi^1\over\partial x_1}(t,\,0,\,0,\,0) &
{\partial \varphi^2\over\partial x_1}(t,\,0,\,0,\,0) &
{\partial \varphi^3\over\partial x_1}(t,\,0,\,0,\,0) \cr
{\partial \varphi^1\over\partial x_2}(t,\,0,\,0,\,0) &
{\partial \varphi^2\over\partial x_2}(t,\,0,\,0,\,0) &
{\partial \varphi^3\over\partial x_2}(t,\,0,\,0,\,0) \cr
{\partial \varphi^1\over\partial x_3}(t,\,0,\,0,\,0) &
{\partial \varphi^2\over\partial x_3}(t,\,0,\,0,\,0) &
{\partial \varphi^3\over\partial x_3}(t,\,0,\,0,\,0) \cr
}
\circ \varphi[t]$$
satisfies
$$\varphi^\ast[t](0, \, 0, \, 0)=(0, \, 0, \, 0),\quad
D\varphi^\ast[t](0,\, 0,\,0) = {\bf 1}_{{\bf R}^3}$$
with
$$|||D^{K}\varphi^\ast[t]|||= |||D^{K}\varphi[t]|||$$
for each $K = 1, \, 2, \, 3$ and each $t$, and
$$\left|\left|\left|
{\partial \varphi^\ast\over \partial t}(t, \, \cdot)\right|\right|\right|
\le 3\bigg(|||D^{0}\varphi||| \cdot 
|||D^{2}\varphi||| + |||D^{1}\varphi[t]|||^2\bigg).$$

\proclaim{{3.6 T}heorem}
There is $C < \infty$ such that the following is true for
all sufficiently small $\delta > 0$. 
Suppose $\gamma_0\co [0, \, \delta] \to \MM$ 
is an arc length parametrization of a length minimizing geodesic 
in $\MM$ and set 
$$\gamma(s, \, t) = \varphi[t]\big(\gamma_0(s)\big)
\quad\hbox{ for each $s$ and $t$}$$
so that $s \to \gamma(s, \, t)$ is an arc length parametrization of
a geodesic in $\MM_t$.  We also set
$$r(s, \, t) = \big|\gamma(0, \, t) - \gamma(s, \, t)\big|
\quad\hbox{ for each $s$ and $t$}$$
and, for (fixed) $0 < R < \delta$, consider
$$ r(R, \, t) = \big|\gamma(0, \, t) - 
\gamma(R, \, t)\big|\quad\hbox{for each $t$.}$$
Then
$${d \over dt}r(R, \,t) = \pm C R^2$$
and
$$\lim_{R \downarrow 0}R^{-1}{d\over dt} r(R, \, t) = 0.$$
\endproc

\prf  We will show
$${d \over dt}r(R, \,t)\bigg|_{t=0} = \pm C R^2.$$

\medskip
{\bf Step 1}\qua
Replacing $\varphi(t, \,p)$ by $\varphi^\ast(t, \,p)
=\varphi(t, \, p) - \varphi(t, \, \gamma_0(0))$ as in 3.4
if necessary we assume without loss of generality that
$\gamma(0, \, t) = (0, \, 0, \, 0)$ for each $t$.

\medskip
{\bf Step 2}\qua  Rotating coordinates if necessary we assume
without loss of generality that ${\bf e}_1$ and ${\bf e}_2$
are tangent to $\MM_0$ at $(0, \, 0, \, 0)$ and that
$\gamma_0'(0) = {\bf e}_1$

\medskip
{\bf Step 3}\qua Rotating coordinates as time changes 
as in 3.5 if necessary
we assume without loss of generality that $D\varphi[t](0, \, 0, \, 0)
= {\bf 1}_{{\bf R}^3}$ for each $t$.

\medskip
{\bf Step 4}\qua We define
$$X(s, \, t) = \gamma(s, \, t) \cdot {\bf e}_1, \quad
Y(s, \, t) = \gamma(s, \, t) \cdot {\bf e}_2, \quad
Z(s, \, t) = \gamma(s, \, t) \cdot {\bf e}_3$$
so that
$$\gamma(s, \, t) = \big(X(s, \, t),\, Y(s, \, t),\, Z(s, \, t)\big)$$
and estimate for each $s$ and $t$:

(a)\qua $X(0, \, t) = Y(0, \, t) = Z(0, \, t) = 0$ (by step 1)

(b)\qua $X_t(0, \, 0) = Y_t(0, \, 0) = Z_t(0, \, 0) = 0$

(c)\qua $X_s(s, \, t)^2 + Y_s(s, \, t)^2 + Z_s(s, \, t)^2 = 1$

(d)\qua $X_s(s, \, t)=\pm 1,\,Y_s(s, \, t) = \pm 1, \, Z_s(s, \, t) = \pm 1$

(e)\qua $1/2 \le r(s, \, t)/|s| \le 1$ (since $\delta$ is small)

(f)\qua $X(s, \, 0) = \pm C s$,  $Y(s, \, 0) = \pm C s$,
$Z(s, \, 0) = \pm C s$

(g)\qua $X_s(0, \, t) = X_s(0, \, 0)$, $Y_s(0, \, t) = Y_s(0, \, 0)$,
$Z_s(0, \, t) = Z_s(0, \, 0)$ (by step 3)

(h)\qua $X_{st}(0, \, 0) = Y_{st}(0, \, 0)= Z_{st}(0, \, 0)=0$

$$X_{st}(s, \, 0) =X_{st}(0, \, 0) + \int_0^sX_{sst}(\eta, \, 0)\, 
d \eta = 0 \pm s \,\sup\big|X_{sst}\big| = \pm Cs, \leqno{(i)}$$
$$Y_{st}(s, \, 0) = \pm Cs, 
\qquad Z_{st}(s, \, 0) = \pm Cs$$

$$X_t(s, \, 0) = X_t(0, \, 0) + \int_0^s X_{st}(\eta, \, 0)\, d\eta
= 0 \pm Cs^2, \leqno{\hbox{(j)}}$$
$$  Y_t(s, \, 0) = \pm Cs^2,\qquad Z_t(s, \, 0) = \pm Cs^2$$ 
(k)\qua $r^2 = X^2 + Y^2 + Z^2$
$$rr_s = XX_s+ YY_s+ ZZ_s, \quad r_s ={1\over r}\big(
XX_s+ YY_s+ ZZ_s \big)\leqno{\hbox{($\ell$)}}$$ 
$$rr_t = XX_t+ YY_t+ ZZ_t,\quad  r_t ={1\over r}\big(
XX_t + YY_t + ZZ_t \big)\leqno{\hbox{(m)}}$$
(n)\qua $r_sr_t + r r_{st} = X_sX_t + X X_{st}+ Y_sY_t + Y Y_{st}+
 Z_sZ_t + Z Z_{st}$
 
(o)\qua evaluating (n) at $t = 0$, $r > 0$ we see
$$\eqalign{{1 \over r(s, \, 0)^2}
\big((\pm Cs)(\pm 1)\big)
\big((\pm Cs)(\pm Cs^2)\big)& + r(s, \, 0) r_{st}(s, \, 0)\cr
=& (\pm 1)(\pm Cs^2) + (\pm Cs)(\pm Cs)\cr}$$
(p)\qua $r_{st}(s, \, 0) = \pm C s$
 
$$r_t(R, \, 0) = r_t (0, \, 0) + \int_0^R r_{st}(s, \, 0)\, ds
=0 +\int_0^R \pm C s\, ds = \pm CR^2.\leqno{\hbox{(q)}}$$
\vglue-.6cm\noindent\hbox{}\hfill$\sq$
%

\proclaim{{3.7 C}orollary}  Suppose triangulation $\TT(j)$ has
maximum edge length $L=L(j)$ and $\langle pq \rangle$ is an
edge in $\TT_1(j)$. Then, for each $t$,
$$\bigg| \varphi[t](p) -\varphi[t](q)\bigg| = \pm C L \quad\hbox{and}\quad
{d \over dt}
\bigg| \varphi[t](p) -\varphi[t](q)\bigg| = \pm C L^2.$$ \endproc

\sh{{3.8\qua S}tabilizing the {f}acets of an {e}dge}
  
Suppose $\TT(j)$ is a triangulation with maximum edge length $L=L(j)$
and that $\langle ABC\rangle, \, \langle ACD\rangle$ are facets
in $\TT_2(j)$ as illustrated
$$\matrix{
&&D=(e,f,0)&& \cr
&\swarrow &&\nwarrow&\cr
(0,0,0)=A&&\longleftrightarrow&&C=(d,0,0)\cr
&\searrow &&\nearrow&\cr
&&B=(a,b,c)&&\cr}.$$
Interchanging $B$ and $D$ if necessary we
assume without loss of generality the the average normal
${\bf n}[0, \,AC]$ to $\MM_0$
at $A$ has positive inner product with $(C-A)\times (D-A)$.
 
\medskip
{\bf1\rm)\qua\bf Fixing $A$ at the origin}\qua Modifying $\varphi$ if necessary
as in 3.4 if necessary we can assume
without loss of generality that $\varphi[t](A) = (0, \, 0, \, 0)$
for each $t$.  As indicated there, various derivative bounds are
increased by, at most, a controlled amount.
 
\medskip
{\bf2\rm)\qua\bf Convenient rotations}\qua  We set
$u(t) = \varphi[t](C), \quad v(t)= \varphi[t](D)$
and use the Gramm--Schmidt orthonormalization process to construct
$$U(t) = {u(t)\over |u(t)|}, \quad 
V(t) = {v(t) - v(t)\cdot U(t)\,U(t)\over
|v(t) - v(t)\cdot U(t)\,U(t)|}, \quad W(t) = U(t) \times V(t).$$
One uses the mean value theorem in checking
$$ |||D^KU(t)|||\le C\left(
\sum_{j=0}^{K+1}|||D^{j}\varphi|||\right), \quad\hbox{etc}$$
for each $K = 0, \, 1, \, 2$.
We denote by $Q(t)$  the orthogonal matrices having columns
equal to $U(t)$, $V(t)$, $W(t)$ respectively
(which is the inverse matrix to its transpose).  Replacing
$\varphi_t$ by $Q(t)\circ \varphi_t$ if necessary, we assume 
without loss of generality that there are functions
$a(t)$, $b(t)$, $c(t)$, $d(t)$, $e(t)$, $f(t)$, 
such that
$$\varphi[t](A) = (0, \, 0, \, 0), \quad
\varphi[t](B) = (a(t), \,b(t), \,c(t)), $$
$$
\varphi[t](C) = (d(t), \, 0, \, 0), \quad
\varphi[t](D) = (e(t), \, f(t), \, 0).$$
We assume without loss of generality the existence of
functions $F[t]\big(x, \, y\big)$ defined for $(x,\,y)$ near $(0, \, 0)$
such that, near $(0, \, 0, \,0)$ our manifold $\MM_t$ is the graph
of $F[t]$.
In particular,
$$c(t) = F[t]\big(a(t), \, b(t)\big).$$
We assert that if $|p| \le CL $, then
$$ |F[t](p)| \le CL^2, \qquad
|\nabla F[t](p)| \le CL.\eqno{(3.8.1)}$$
To see this, first we note that $F[t](A) = F[t](C) = F[t](D) = 0$.
Next we invoke Rolle's theorem to conclude the
existence of $c_1$ on segment $AD$ and $c_2$ on segment
$CD$ such
$$\left<{D-A\over |D-A|}, \, DF[t](c_1)\right> = 0 
=\left<{D-C\over |D-C|}, \, DF[t](c_2)\right> .$$
Since $|p| \le CL $ we infer
$$\left<{D-A\over |D-A|}, \, DF[t](p)\right> =  \pm CL, \qquad
\left<{D-C\over |D-C|}, \, DF[t](p)\right> = \pm CL .$$
In view of 2.1.6(vi)(vii)(viii) and 2.2.7 we infer that 
${\bf e}_1$ and ${\bf e}_2$ 
are bounded linear combinations of
$(D-A)/|D-A|$ and $(D-C)/ |D-C|$ from which we conclude that
$|\nabla F[t](p)| \le CL$.  This in turn implies that
$ |F[t](p)| \le CL^2$ as asserted.

Since
$${\partial \over \partial t}F[t](0, \, 0) = 0$$
we infer
$${\partial \over \partial t}F[t](p)= \pm  C L\eqno{(3.8.2)}$$
and since
$${\partial \over \partial t}\left(\varphi[t](A)\cdot
 {\bf e}_3\right) = 0$$
we infer
$$c'(t) = {\partial \over \partial t}F[t](a(t), \, b(t))
={\partial \over \partial t}\left(\varphi[t](B)\cdot {\bf e}_3\right)
= \pm  C L.\eqno{(3.8.3)}$$

\proclaim{{3.9 P}roposition}
Let $L, \, A, \, B, \, C, \, D, \, a, \, b,\, c,\, d, \, e, \, f$ be
as in 3.8.  Then

{\rm(1)}\qua $a'(t) = \pm CL^2$

{\rm(2)}\qua $b'(t) = \pm CL^2$

{\rm(3)}\qua $c'(t) = \pm CL$

{\rm(4)}\qua $d'(t) = \pm CL^2$

{\rm(5)}\qua $e'(t) = \pm CL^2$

{(\rm6)}\qua $f'(t) = \pm CL^2$.\endproc

\prf  According to 3.7, if
$r(t)$ denotes the distance between the endpoints of
an edge of arc length $L$ at time $t$, then
$$r'(t) = \pm C L^2.$$

(i)\qua We invoke 3.7 directly to infer (4) above.

(ii)\qua We apply 3.7 to the distance between
$(0,\, 0, \, 0)$ and $(e, \, f,\, 0)$ to infer
$${d \over dt}\big(e^2 + f^2\big)^{1\over 2}
={ \big(e e' + ff'\big) 
 \over \big(e^2 + f^2\big)^{1\over 2}}
= \pm CL^2,
\qquad
ee'+ ff' = \pm C L^3.$$
(iii)\qua We apply 3.7 to the distance between
$(d,\, 0, \, 0)$ and $(e, \, f,\, 0)$ to infer
$$\eqalign{{d \over dt}\big((e-d)^2 + f^2\big)^{1\over 2}
=&{\big(e-d)(e'-d') + ff'\big)
 \over \big((e-d)^2 + f^2\big)^{1\over 2}}
 = \pm CL^2,\cr
(e -& d)(e'- d')+ ff' = \pm C L^3.\cr}$$
We subtract the first inequality from the second to infer
$$ed' - d e' + d d' = \pm CL^3, \qquad de' \pm CL^3, \qquad e'=\pm CL^2.$$
Assertions (5) and (6) follow readily.

(iv)\qua We apply 3.7 to the distance between
$(0,\, 0, \, 0)$ and $(a, \, b,\, c)$ to infer
$${d \over dt}\big(a^2 + b^2 + c^2\big)^{1\over 2}
={\big(a a' + bb' + c c'\big)
\over\big(a^2 + b^2 + c^2\big)^{1\over 2}}
 = \pm CL^2,
\qquad
aa'+ bb'+ cc' = \pm C L^3.$$
(v)\qua We apply 3.7 to the distance between
$(d,\, 0, \, 0)$ and $(a, \, b,\, c)$ to infer
$${d \over dt}\big((a-d)^2 + b^2 + c^2\big)^{1\over 2}
={\big((a-d)(a'-d') + bb' + c c'\big)
 \over \big((a-d)^2 + b^2+c^2\big)^{1\over 2}}
 = \pm CL^2,$$
$$ (a - d)(a'- d')+ bb' + c c' = \pm C L^3.$$
We subtract the first inequality form the second to infer
$$ad' - d a' + d d' = \pm CL^3, \qquad da' \pm CL^3, \qquad a'=\pm CL^2,$$
which gives assertion (1).

(vi)\qua We estimate from 3.8 that
$$c = F[t](a, \, b) = \pm CL^2, \quad
c'= {d \over dt}F[t](a, \, b) + \nabla F[t](a, \, b) \cdot (a', \, b')
=\pm CL,$$
which gives (3) above.  We have also
$cc'= \pm C L^3$.
We recall (iv) above and estimate
$$a a' +bb' + cc' = \pm C L^3, \quad bb' =\pm CL^3, \quad b'= \pm CL^2,$$
which is (2) above.
\E

\proclaim{{3.10 P}roposition}  Suppose $\TT(j)$ is a triangulation with
maximum edge length $L=L(j)$ and
$\langle pq\rangle$ is an edge in $\TT_1(j)$.
Abbreviate $\theta(t) = \theta[t, \, j](pq)$.  Then, for each $t$,
$$\theta(t) = \pm CL\leqno{\rm(1)}$$
$$2 \sin\left({\theta(t) \over 2}\right) = \pm CL\leqno{\rm(2)}$$
$$\theta'(t) = \pm C\leqno{\rm(3)}$$
$${d \over dt}\left[2 \sin\left({\theta(t) \over 2}\right)\right]
 = \pm C\leqno{(4)}$$
$${d \over dt}\left[2 \sin\left({\theta(t) \over 2}\right)- \theta\right]
 = \pm CL^2.\leqno{(5)}$$\endproc

\prf 
Making the modifications of 3.8 if necessary, we assume without
loss of generality (in the terminology there) that
$\varphi[t](p) = A = (0, \, 0, \,0)$, 
$\varphi[t](q) = C = (d(t), \, 0, \,0)$, 
and that there are $\langle pqB_\ast\rangle, \, \langle pqD_\ast\rangle
\in \TT_2(j)_0$ with
$\varphi[t](B_\ast) = B = (a(t), \, b(t), \,c(t))$, 
$\varphi[t](D_\ast) = D = (e(t), \, f(t), \,0)$. 

The unit normal to $\overline{ACD}$ is $(0, \, 0, \, 1)$ while the 
unit normal to $\overline{ABC}$ is
$${(0, \,-c, \, b)\over (b^2 + c^2)^{1\over 2}}$$
so that\qquad
$\displaystyle\cos\theta = {b \over (b^2 + c^2)^{1\over 2}}$, 
$$
\sin\theta =\pm \left(1 - \cos^2\theta\right)^{1\over 2} 
= \pm\left( 1 - {b^2 \over b^2 + c^2}\right)^{1\over 2}
= \pm
{c \over (b^2 + c^2)^{1\over 2}}= \pm CL$$
in view of 3.8.  Assertions (1) and (2) follow.
We compute further
$$(\sin\theta)' = \cos \theta\, \theta' =
\pm{(b^2+ c^2)^{1\over 2}c' - c{bb' + cc'\over (b^2 + c^2)^{1\over 2}}
\over b^2 + c^2} = \pm C$$
in view of 3.9(1)(2)(3) and 3.8.  Assertion (3) and (4) follow.
Assertion (5) follows from differentiation and assertions (1) and (3).
\E

\proclaim{{3.11  P}roposition}  Suppose $\TT(j)$ is a triangulation with
maximum edge length $L=L(j)$ and  $\langle pq \rangle$ is an edge in
$\TT_1(j)$.  Then

{\rm(1)}\qua  ${\bf n}[t,\, j](pq) = \big(0, \, \pm CL, \, 1 \pm C L^4\big)$ 

{\rm(2)}\qua  $(d / dt) \big({\bf n}[t,\, j](pq)\big)
 = \big(0, \, \pm C, \, \pm CL\big) +  \big(\pm CL, \, 
\pm CL, \, \pm CL\big)$

{\rm(3)}\qua  $g[t, \, j](pq) = \big( \pm CL, \, \pm CL, \,1 \pm CL^2\big)$

{\rm(4)}\qua  $(d/dt) g[t, \, j](pq) = \big(\pm C, \, \pm C, \, 0\big) + \big(
\pm CL, \, \pm CL, \, \pm CL\big)$ 

{\rm(5)}\qua  ${\bf n}[t,\, j](pq)\cdot g[t, \, j](pq) = 1 \pm CL^2$

{\rm(6)}\qua  $(d / dt)\bigg(
{\bf n}[t,\, j](pq)\cdot g[t, \, j](pq)\bigg)
=\pm CL$

{\rm(7)}\qua  $ 1 -{\bf n}[t,\, j](pq) \cdot g[t, \, j](pq) = \pm CL^2$.
\endproc

\prf  We let $A$, $B$, $C$, $D$, $F[t]$, $b(t)$, $c(t)$, $d(t)$
be as in 3.8.
We abbreviate ${\bf n} ={\bf n}[t, \, j](pq)$ and estimate
$$\eqalign{{\bf n}&= 
{(0, \, 0, \, 1) + (0, \, -c, \, b)/(b^2 + c^2)^{1\over 2}
\over
\big| (0, \, 0, \, 1) + (0, \, -c, \, b)/(b^2 + c^2)^{1\over 2} \big|}\cr
&={\big(0, \, -c, \, b + (b^2 + c^2)^{1\over 2}\big) \over
2^{1\over 2}\big(b^2 + c^2 + b(b^2 + c^2)^{1\over 2}
\big)^{1\over 2}}.\cr}$$
The first assertion follows from 3.8.1.
We differentiate to conclude\qua ${\bf n}'=$
$$\eqalign{&{\pm CL \big(0, \, -c', \, b' \pm C(b b' + c c')/L
- (L/L)\big(bb' + cc' \pm b'L + \pm C(b/L)(bb' + cc')\over
\pm  L^2}\cr
&\phantom{XX}= \big(0, \, \pm C, \, \pm CL\big) + \big(\pm CL, \, \pm CL, \, 
\pm CL\big)\cr}$$
in view of 3.9(2)(3).  This is assertion (2).

We abbreviate $g = g[t, \, j](pq)$ and estimate
$$\eqalign{g &= {1\over d(t)}\int_0^{d(t)}{
\big(-F[t]_x, \, -F[t]_y, \, 1\big)\over
\big|\big(-F[t]_x,\, -F[t]_y,\, 1\big)\big|}\cr
&= {1\over d(t)}\int_0^{d(t)}{
\big(-F[t]_x, \, -F[t]_y, \, 1\big)\over
\bigg(\big(F[t]_x^2 F[t]_y^2+1\bigg)^{1\over 2}}.\cr}$$
The third assertion follows from 3.8.1.  We differentiate to estimate
that ${d g /dt}$ equals
$$\eqalign{
 &{-d'\over d^2} \int_0^{d(t)}{
\big(-F[t]_x, \, -F[t]_y, \, 1\big)\over
\big(1 +F[t]_x^2+ F[t]_y^2\big)^{1\over 2}}
+ {d'\over d}
{\big(-F[t]_x, \, -F[t]_y, \, 1\big)\over
\big(1 + F[t]_x^2+ F[t]_y^2\big)^{1\over 2}}\cr 
&
\phantom{XXX}+
{1\over d}\int_0^d{\pm CL\big(-F[t]_{tx}, \, -F[t]_{ty}, \, 0\big)
\over 1 + F[t]_x^2 + F[t]_y^2}\cr
&\phantom{XXX}
-{1\over d}\int_0^d {\big(-F[t]_x,\, -F[t]_y, \, 1\big)
(\pm C/L)\big(F[t]_xF[t]_{tx}
+F[t]_yF[t]_{ty}\big)
\over 1 + F[t]_x^2 + F[t]_y^2}=\cr
 &L\big(\pm C, \, \pm C, \, \pm C\big)
+ L\big(\pm C, \, \pm C, \, \pm C\big)
+ \big(\pm C, \, \pm C, \, 0\big)
+ L\big(\pm C, \, \pm C, \, \pm C\big)\cr
}$$
which gives assertion (4).  Assertion (5) follows from assertions (1) and
(3). Assertion (6) follows from assertions
(1), (2), (3), (4) and integration by parts.
Assertion (7) follows from assertions (1) and (3).  \E

\section{Constancy of the mean curvature integral}

\sh{{4.1\qua  T}he {d}erivative {e}stimates}

Suppose triangulation $\TT(j)$ has maximum edge length $L= L(j)$.
We recall from 2.2.10 that
$$\eqalign{&\delta V[\TT(j, \, t)]\big(g[t]\big)\cr 
&= \sum_{\langle pq \rangle \in \TT_1(j)}
\bigg[ \big|\varphi[t](p) - \varphi[t](q)\big|
\bigg]\,\bigg[2\, \sin 
\left({\theta[t,\,j](pq)\over 2}\right)\bigg]
\, \bigg[{\bf n}[t, \, j](pq)\cdot g[t, \, j](pq)\bigg]\cr}$$
and we estimate, for each $t$ that
$$\eqalign{&{d \over d t}\bigg(\delta V[\TT(j)_t]\big(g[t]\big)\bigg)
\cr & 
= \sum_{\langle pq \rangle \in \TT_1(j)}
\bigg[ \big|\varphi[t](p) - \varphi[t](q)\big| \bigg]'\,\bigg[2\, \sin 
\left({\theta[t,\,j](pq)\over 2}\right)\bigg]
\, \bigg[{\bf n}[t, \, j](pq)\cdot g[t, \, j](pq)\bigg]\cr
& 
+ \sum_{\langle pq \rangle \in \TT_1(j)}
\bigg[ \big|\varphi[t](p) - \varphi[t](q)\big|\bigg]
\,\bigg[2\, \sin 
\left({\theta[t,\,j](pq)\over 2}\right)\bigg]'
\, \bigg[{\bf n}[t, \, j](pq)\cdot g[t, \, j](pq)\bigg] \cr
& 
+ \sum_{\langle pq \rangle \in \TT_1(j)}
\bigg[ \big|\varphi[t](p) - \varphi[t](q)\big|\bigg]
\,\bigg[2\, \sin 
\left({\theta[t,\,j](pq)\over 2}\right)\bigg]
\, \bigg[{\bf n}[t, \, j](pq)\cdot g[t, \, j](pq)\bigg]'.\cr}$$

We assert that
$${d \over dt}\bigg(\delta V[\TT(j, \, t)]\big(g[t]\big) \bigg)
= \sum_{\langle pq \rangle \in \TT_1(j)}\pm CL^3
= \sum_{\langle pq \rangle \in \TT_1(j)}\pm CL(j)^3.$$
To see this we will estimate each of the three summands above.

{\bf First summand}\qua We use 3.7, 3.10(2), 3.11(5) to estimate for each
 $pq$,
$$\eqalign{\bigg[ \big|\varphi[t](p) - \varphi[t](q)\big| \bigg]'\,\bigg[2\, \sin 
&\left({\theta[t,\,j](pq)\over 2}\right)\bigg]
\, \bigg[{\bf n}[t, \, j](pq)\cdot g[t, \, j](pq)\bigg]
\cr 
=&\big(CL^2\big)\big(CL\big)\big(1 \pm CL^2\big).\cr }$$

{\bf Second summand}\qua We use 3.10(5), 3.11(7) to estimate for each $pq$,
$$\eqalign{\bigg[& \big|\varphi[t](p) - \varphi[t](q)\big|\bigg]
\,\bigg[2\, \sin \left({\theta[t,\,j](pq)\over 2}\right)\bigg]'
\, \bigg[{\bf n}[t, \, j](pq)\cdot g[t, \, j](pq)\bigg] \cr
&= \bigg[ \big|\varphi[t](p) - \varphi[t](q)\big|\bigg]
\bigg[ \theta[t,\,j](pq)\bigg]'\cr
&\phantom{X}
+ \bigg[ \big|\varphi[t](p) - \varphi[t](q)\big|\bigg]
\,\bigg[2\, \sin \left({\theta[t,\,j](pq)\over 2}\right)
-\theta[t,\,j](pq)\bigg]'\cr
&\phantom{X}
+\bigg[ \big|\varphi[t](p) - \varphi[t](q)\big|\bigg]
\,\bigg[2\, \sin \left({\theta[t,\,j](pq)\over 2}\right)\bigg]'
\, \bigg[{\bf n}[t, \, j](pq)\cdot g[t, \, j](pq) - 1\bigg]\cr
&= \bigg[ \big|\varphi[t](p) - \varphi[t](q)\big|\bigg]
\bigg[ \theta[t,\,j](pq)\bigg]'
\pm \big(CL\big)\big(CL^2\big)\pm 
\big(CL\big)\big(C\big)\big(CL^2\big).
\cr}$$
 
{\bf Third summand}\qua  We use 3.10(2) and 3.11(6) to estimate
$$\eqalign{\bigg[ \big|\varphi[t](p) - \varphi[t](q)\big|\bigg]
\,\bigg[2\,& \sin \left({\theta[t,\,j](pq)\over 2}\right)\bigg]
\, \bigg[{\bf n}[t, \, j](pq)\cdot g[t, \, j](pq)\bigg]'\cr
&
=\big(CL\big) \big(CL\big) \big(CL\big).\cr}$$
According to Schlafli's formula [\S],
$$\sum_{\langle pq \rangle \in \TT_1(j)}
\bigg[ \big|\varphi[t](p) - \varphi[t](q)\big|\bigg]
\bigg[ \theta[t,\,j](pq)\bigg]'= 0.$$
Our assertion follows.

\proclaim{{4.2 M}ain {T}heorem}

{\rm(1)}\qua  For each fixed time $t$, 
$$\lim_{j \to \infty}\delta V[\TT(j, \, t)]\big(g[t]\big) = \delta
V_t\big(g[t]\big).$$

{\rm(2)}\qua  For each fixed $j$, 
$\delta V[\TT(j)_t]\big(g[t]\big)$ is a differentiable function of $t$
and 
$$\lim_{j \to \infty}{d \over dt}
\bigg(\delta V[\TT(j)_t]\big(g[t]\big)\bigg)=0$$
uniformly in $t$.

{\rm(3)}\qua  For each $t$
$$\int_{\MM_t} H_t \, d {\cal H}^2 =
\int_{\MM} H \, d {\cal H}^2.$$
This is the main result of this note.

\prf
To prove the first assertion, we check that 
$$(\rho_t)_\sharp V[\TT(j, \,t)] = V_t$$
for each $t$ and all large $j$.  Indeed, the $\tau$ regularity of our
triangulations implies that the normal directions of the $\NN[\TT(j)_t]$
are very nearly equal to the normal directions of nearby points on
$\MM_t$ and that the restriction of $D\rho_t$ to the tangent planes
of the $\NN[\TT(j)_t]$ is very nearly an orthogonal injection.  The
first assertion follows with use of the first variation formula given
in [\AW 4.1, 4.2].
Assertion (2) follows from 4.1 since
$$ \sum_{\langle pq \rangle \in \TT_1(j)} L(j)^2$$
is dominated by the area of $\MM$ (see 2.2.12) and $ \lim_{j \to
\infty}L(j) = 0$.
Assertion (3) follows from assertions (1) and (2) and our
observation in 2.1.4. \E  

\rk{Acknowledgements}
Fred Almgren tragically passed away shortly after this note was
written. Since then, the main result for smooth surfaces has been
reproved in an easier way and generalized to the setting of Einstein
manifolds by J-M Schlenker together with the second author of the
current paper [\RS]. Nonetheless, it seems clear that the methods
used here can be used to extend these results in other directions.

\references\Addresses\recd\bye